\documentclass[12pt]{article}
\usepackage{amssymb}
\usepackage{a4}
\begin{document}
\def\st{\, : \,}
\def\kbar{{\mathchar'26\mkern-9muk}}  
\def\bra#1{\langle #1 \vert}
\def\ket#1{\vert #1 \rangle}
\def\vev#1{\langle #1 \rangle}
\def\ad{\mbox{ad}\,}
\def\ker{\mbox{Ker}\,}
\def\im{\mbox{Im}\,}
\def\der{\mbox{Der}\,}
\def\ad{\mbox{ad}\,}
\def\b#1{{\mathbb #1}}
\def\c#1{{\cal #1}}
\def\pt{\partial_t}
\def\px{\partial_1}
\def\bpx{\bar\partial_1}
\def\la{\langle}
\def\ra{\rangle}
\def\nn{\nonumber \\}
\def\pn{\par\noindent}
\def\etal{{\it et al.}\ }
\def\sq{\mbox{\rlap{$\sqcap$}$\sqcup$}}
\def\R{{\cal R}\,}
\newcommand{\tr}{\triangleright\,}
\newcommand{\tl}{\,\triangleleft}
\newcommand{\tro}{\triangleright^{op}\,}
\newcommand{\tlo}{\,\stackrel{op}{\triangleleft}}
\def\cross{{\triangleright\!\!\!<}}
\def\cocross{{>\!\!\!\triangleleft\,}}
\def\uqg{\mbox{$U_q{\/\bf g}$ }}
\def\uqso{\mbox{$U^{op}_q{\/so(N)}$ }}
\def\uqsp{\mbox{$U_q^+{\/so(N)}$ }}
\def\uqsn{\mbox{$U_q^-{\/so(N)}$ }}
\def\uqs{\mbox{$U_q{\/so(N)}$ }}

\renewcommand{\thefootnote}{\fnsymbol{footnote}}

\renewcommand{\theequation}{\thesection.\arabic{equation}}
\newcommand{\initiate}{\setcounter{equation}{0}}

\newcommand{\ba}{\begin{array}}
\newcommand{\ea}{\end{array}}
\newcommand{\be}{\begin{equation}}
\newcommand{\ee}{\end{equation}}
\newcommand{\bea}{\begin{eqnarray}}
\newcommand{\eea}{\end{eqnarray}}
\newcommand{\beas}{\begin{eqnarray*}}
\newcommand{\eeas}{\end{eqnarray*}}
%
%
%
\newtheorem{prop}{Proposition}
\newtheorem{lemma}{Lemma}
\newtheorem{theorem}{Theorem}
\newtheorem{corollary}{Corollary}
%
%
\newenvironment{proof}[1]{\vspace{5pt}\noindent{\bf Proof #1}\hspace{6pt}}%
{\hfill\sq}
\newcommand{\bp}{\begin{proof}}
\newcommand{\ep}{\end{proof}\par\vspace{10pt}\noindent}
%
%

\title{Frame formalism for the N-dimensional quantum Euclidean spaces}

\author{B. L. Cerchiai,$\strut^{1,2}$ \,
        G. Fiore,$\strut^{3,4}$ \, J. Madore$\strut^{5,2}$ \\\\
        \and
        $\strut^1$Sektion Physik, Ludwig-Maximilians-Universit\"at,\\
        Theresienstra\ss e 37, D-80333 M\"unchen
        \and
        $\strut^2$Max-Planck-Institut f\"ur Physik\\
        F\"ohringer Ring 6, D-80805 M\"unchen
        \and
        $\strut^3$Dip. di Matematica e Applicazioni, Fac.  di Ingegneria\\ 
        Universit\`a di Napoli, V. Claudio 21, 80125 Napoli
        \and
        $\strut^4$I.N.F.N., Sezione di Napoli,\\
        Mostra d'Oltremare, Pad. 19, 80125 Napoli
        \and
        $\strut^5$Laboratoire de Physique Th\'eorique et Hautes Energies\\
        Universit\'e de Paris-Sud, B\^atiment 211, F-91405 Orsay
        }
\date{}

\maketitle
\abstract{We sketch our application~\cite{CerFioMad00} of a 
non-commutative version of the Cartan `moving-frame' formalism to the quantum
Euclidean space $\b R^N_q$, the space which is covariant under the action
of the quantum group $SO_q(N)$. For each of the two covariant differential
calculi over $\b R^N_q$ based on the $R$-matrix formalism, we summarize
our construction of
a frame, the dual inner derivations, a metric and two torsion-free almost
metric compatible covariant derivatives with a vanishing curvature.
To obtain these results we have developed a technique which fully exploits the 
quantum group covariance of $\b R^N_q$. We first find a frame in the larger
algebra $\Omega^*(\b R^N_q) \cocross \uqs$. Then we define homomorphisms
from $\b R^N_q \cocross U_q^{\pm}{\/so(N)}$ to $\b R^N_q$  which we use to
project this frame in $\Omega^*(\b R^N_q)$.
}

\vfill
\noindent
Preprint 00-29 Dip. Matematica e Applicazioni, Universit\`a di Napoli
\newpage

\section{Introduction and preliminaries}

Non-commutative geometry as a way to describe the structure of space-time 
at small distances has been proposed since 1947~\cite{Sny47}. 
The original claim however made by him that it would serve as a
Lorentz-invariant cut-off has been recently the object of some
controversy.  Here we briefly
describe our application~\cite{CerFioMad00} of a non-commutative
generalization~\cite{DimMad96} of the Cartan `moving-frame' formalism to the
quantum Euclidean spaces $\b R^N_q$~\cite{FadResTak89}, the spaces which are 
covariant under the quantum group $SO_q(N)$.

In Section~\ref{preli1} we recall the basic concepts of 
non-commutative geometry,
i.e. given a non-commutative algebra and a differential
calculus on it how one can
introduce the notions of a moving-frame or `Stehbein'~\cite{DimMad96}, of
a corresponding metric and of a covariant derivative. For a more detailed 
exposition see~\cite{Mad95}.
In Section~\ref{preli2} we review the definition of the quantum Euclidean
spaces~\cite{FadResTak89} and of the differential calculi~\cite{CarSchWat91,
Ogi92, WesZum90} on them. 
Finally in Section~\ref{applic} we show how to construct the frame on $\b R^N_q$ and
with its help the corresponding metric and covariant derivatives. 
It is necessary
to enlarge the algebra with a `dilatator', the square roots and
the inverses of some
elements and in the case of even $N$ also with one of the components of the
angular momentum.

The quantum group covariance is an essential ingredient of our construction.
We first define a frame in the cross-product $\Omega^1(\b R^N_q) \cocross 
\uqs$ and then project it to a frame in $\Omega^1(\b R^N_q)$ through
homomorphisms $\varphi^{\pm}$ from $\b R^N_q \cocross U_q^{\pm}{\/so(N)}$ to
$\b R^N_q$,  which we apply to the components.
In other words the components of the frame
in the $dx^i$ basis automatically provide a `local realization' of 
$U_q^{\pm}{\/so(N)}$ in the extended algebra of $\b{R}^N_q$, i.e.. they satisfy 
the `RLL' and the `gLL' relations fulfilled by the 
Faddeev-Reshetikhin-Takhtadjan~\cite{FadResTak89} 
generators of $U_q^{\pm}{\/so(N)}$.
For odd $N$ we have the interesting result that $\varphi^{\pm}$ can be glued
to a homomorphism $\varphi:\b R^N_q \cocross U_q{\/so(N)} \rightarrow R^N_q$. 

We recover for $\b R^N_q$ the formal `Dirac operator'~\cite{Con94}, as it
had already been found in~\cite{Zum97, Ste96}. In this way we construct
a link between the approach to noncommutative geometry
of Woronowicz~\cite{Wor}, which is based on the
quantum group covariance, and the one of Connes~\cite{Con94}, which is based
on the `Dirac operator'. Moreover, this method possibly suggests the correct 
choice of the physical coordinates~\cite{FioMad00}. 
For an application of this formalism to the development of a quantum
field theory, see e.g. \cite{GroMadSte00}. Note that 
the role of frame in the Woronowicz 
bicovariant differential calculi is played by the left- (or right-) invariant
1-forms; this has been applied to multiparametric deformations (including
as a particular case the one-parameter one at hand) of the quantum 
Euclidean space in Ref. \cite{AscCast96,AscCastSca97}.

\subsection{The Cartan moving-frame formalism}
\label{preli1}

Following \cite{DimMad96}, we briefly review here a non-commutative 
generalization of the moving-frame formalism of E. Cartan. The starting 
points are a  noncommutative algebra $\c{A}$, which in the commutative 
limit should reduce to the algebra of functions on a parallelizable
manifold $M$, and a differential calculus $\Omega^*(\c{A})$ on it,
which in the same limit should become the de~Rham differential calculus on
$M$.  The set of $1$-forms $\Omega^1(\c{A})$ is assumed then to be a free
$\c{A}$-module of rank $N$. In addition we postulate the existence of a
particular basis, a `frame'
or `Stehbein', which we denote with $\{ \theta^a \}_{1 \le a \le N}$ and  
which has the property of commuting with the elements of $\c{A}$:
\be
[f,\theta^a]=0.                                              \label{thetaf}
\ee
This basis $\theta^a$ is also required to be dual to a set of inner
derivations $e_a=\ad \lambda_a$:
\be
df=e_a f \theta^a=[\lambda_a,f] \theta^a \:\:\: \forall f \in \c{A}. \label{df}
\ee
The integer $N$ can be interpreted as the dimension of $M$.

Under these assumptions a formal `Dirac operator'~\cite{Con94} can be defined
\be
\theta = -\lambda_a \theta^a.                                  \label{dirac}
\ee
It has the property
\be
df = -[\theta,f].
\ee
We shall require the center $\c{Z}(\c{A})$ of $\c{A}$ to be trivial:
$\c{Z}(\c{A})=\b{C}$. If necessary, one can enlarge
the algebra until it is so. 

Using the frame, the relations constraining the wedge product
(we omit the symbol $\wedge$) will take the form
\be
\theta^a \theta^b=P^{ab}{}_{cd} \theta^c \theta^d ,
\qquad P^{ab}{}_{cd} \in \c{Z}(\c{A}).\label{prod}    
\ee 
The matrix $P$ is necessarily a projector, and goes to the
antisymmetric one in the commutative limit. We shall denote by 
$\pi:\Omega^1(\c{A})\otimes\Omega^1(\c{A})\to\Omega^2(\c{A})$ 
the $\c{A}$-bilinear projector such that 
$\pi(\theta^a\otimes \theta^b)=\theta^a \theta^b$.
Then for the condition $d^2 = 0$ to hold 
the $\lambda_a$ have to satisfy a quadratic relation of the form
\be
2 \lambda_c \lambda_d P^{cd}{}_{ab} - 
\lambda_c F^c{}_{ab} - K_{ab} = 0   ,
\qquad F^c{}_{ab}, \: K_{ab}  \in \c{Z}(\c{A}).          \label{manca}
\ee
For the quantum Euclidean spaces $\b{R}^N_q$ the linear and constant terms
will vanish:
\be
F^c{}_{ab}=K_{ab}=0.                          \label{linear}
\ee
The metric is defined as a non-degenerate $\c{A}$-bilinear map
\be
g: \Omega^1(\c{A}) \otimes_{\c{A}} \Omega^1(\c{A}) \rightarrow \c{A}. 
         \label{metric}
\ee
It is completely determined if we assign its value in a basis
\be
g(\theta^a \otimes \theta^b)=:g^{ab}.
\ee
As a further step, one has to introduce a `deformed flip', i.e. an
$\c{A}$-bilinear map \cite{Mou95}
\be
\sigma: \Omega^1(\c{A}) \otimes_\c{A} \Omega^1(\c{A}) \rightarrow
\Omega^1(\c{A}) \otimes_\c{A} \Omega^1(\c{A}),
\qquad \sigma (\theta^a \otimes \theta^b) = :
S^{ab}{}_{cd} \theta^c \otimes \theta^d                 \label{2.2.4}
\ee
which we suppose~\cite{FioMad98} to satisfy the braid relation.
Bilinearity implies that $g^{ab} \in \c{Z}(\c{A})=\b{C}$ and
$S^{ab}{}_{cd}\in \c{Z}(\c{A})=\b{C}$.
Notice that in the classical limit the bilinearity condition, which 
e.g. for the metric reads, 
\be
g(f\xi\otimes\eta h)=f\: g(\xi\otimes\eta)\: h,
\ee
amounts to locality
\[
[g(f\xi\otimes\eta h)](x)=f(x)\: [g(\xi\otimes\eta)](x)\: h(x),
\]
what can be seen as a rationale for the definition given above.
The flip is necessary in order to construct a covariant derivative
$D$~\cite{DubMadMasMou96}, i.e.. a map
\be
D: \Omega^1(\c{A}) \rightarrow \Omega^1(\c{A}) \otimes \Omega^1(\c{A})
\ee
satisfying a left and right Leibniz rule:
\be
D (f \xi) =  df \otimes \xi + f D\xi, \qquad
D(\xi f) = \sigma (\xi \otimes df) + (D\xi) f. 
\ee
Then a consistent torsion map can be found
\be
\Theta:\Omega^1(\c{A})\rightarrow\Omega^2(\c{A}), \qquad \Theta=d-\pi\circ D.
\ee
where bilinearity requires that
\be
\pi \circ (\sigma + 1) = 0.               \label{torsion}
\ee
The associated curvature map is then given by
\be
\mbox{Curv} \equiv D^2= \pi_{12} \circ D_2 \circ D, \qquad
\mbox{Curv} (\theta^a) =
- {1 \over 2} R^a{}_{bcd} \theta^c \theta^d \otimes \theta^b.    \label{curv}
\ee
Here $D_2$ is a continuation of the map (\ref{2.2.4}) to the 
tensor product $\Omega^1(\c{A}) \otimes_\c{A} \Omega^1(\c{A})$
\be
D_2(\xi \otimes \eta) = D\xi \otimes \eta + \sigma_{12} (\xi \otimes D\eta).
\ee
In the case (\ref{linear}), a torsion-free covariant
derivative~\cite{DubMadMasMou96} is
\be
D \xi = - \theta \otimes \xi +\sigma (\xi \otimes \theta).  \label{covdev0}  
\ee
The covariant derivative is compatible with the metric
if $g_{23}\circ D_2= d\circ g$. Under the assumption (\ref{linear})
this is equivalent to~\cite{DubMadMasMou95}
\be
S^{ae}{}_{df} g^{fg} S^{bc}{}_{eg} = g^{ab} \delta^c_d. \label{met-comp}
\ee

\subsection{The quantum Euclidean spaces} \label{preliminaries}
\label{preli2}

Now, following~\cite{FadResTak89}, we briefly review 
the definition and some of the main results about the $N$-dimensional
quantum Euclidean space $\b{R}^N_q$.

The first building block is the matrix $\hat{R}$~\cite{FadResTak89}
for $SO_q(N,\b{C})$ and its projector decomposition:
\be
\hat R = q\c{P}_s - q^{-1}\c{P}_a + q^{1-N}\c{P}_t.       \label{projectorR}
\ee
with $\c{P}_s$, $\c{P}_a$, $\c{P}_t$ $SO_q(N)$-covariant
$q$-deformations of the symmetric trace-free,
antisymmetric and trace projectors respectively. 
The matrix $\c{P}_t$ can be expressed using the metric matrix
$g_{ij}$. Explicitly the latter reads
\be
g_{ij}=g^{ij}=q^{-\rho_i} \delta_{i,-j}.
\ee
It is a $SO_q(N)$-isotropic tensor and is a deformation of the 
ordinary Euclidean metric.
Here and in the sequel $n$ is the rank of $SO(N,\b{C})$, the indices take the 
values
$i=-n,\ldots,-1,0,1,\ldots n$ for $N$ odd,
and $i=-n,\ldots,-1, 1,\ldots n$ for $N$ even.
Moreover, we have introduced the notation
$(\rho_i)=(n-\frac{1}{2},\ldots,\frac{1}{2},0,-\frac{1}{2},\ldots,\frac{1}{2}-n)$
for $N$ odd, $(n-1,\ldots,0,0,\ldots,1-n)$ for $N$ even. 

The $N$-dimensional quantum Euclidean space is defined as the associative
algebra $\b{R}^N_q$ generated by elements $\{x^i\}_{i=-n,\ldots,n}$ with
relations
\be
\c{P}_{a}{}^{ij}_{kl} x^k x^l=0.                              \label{xrel}
\ee
This choice for the values of the indices is the most natural, 
because it is directly related to the weight vector associated
to the Cartan subalgebra of the quantum group.

There are~\cite{CarSchWat91} two $SO_q(N)$-co\-vari\-ant quadratic 
differential calculi on $\b R^N_q$: $\Omega^*(\b{R}_q^N)$ and 
$\bar \Omega^*(\b{R}_q^N)$. They are generated by $\xi^i=dx^i$ and 
$\bar{\xi}^i=\bar{d} x^i$ respectively, with relations
\bea
x^i \xi^j = q\,\hat R^{ij}_{kl} \xi^k x^l,\quad &\quad \c{P}_a{}_{kl}^{ij}
\xi^k \xi^l = \xi^i\xi^j
&\:\:\mbox{for } \Omega^1(\b{R}_q^N),               \label{xxirel}\\
x^i \bar \xi^j = q^{-1}\,\hat R^{-1}{}^{ij}_{kl} \bar\xi^k x^l,    
& \quad \c{P}_a{}_{kl}^{ij} \bar \xi^k \bar \xi^l = \xi^i\xi^j
&\:\:\mbox{for } \bar{\Omega}^1(\b{R}_q^N).           \label{xxistar}
\eea
The nilpotency of $d$ then automatically fixes the wedge product 
relations respectively among the $\xi^i$ and the $\bar\xi^i$.

For $q \in \b{R}^+$ $\b{R}_q^N$ can be equipped with a consistent conjugation
with the help of the metric to lower and raise the indices 
$(x^i)^*= x^j g_{ji}$.
In this way we obtain what is known as real quantum Euclidean space.
The $*$-structure can be extended to $\Omega^1(\b{R}_q^N) \oplus
\bar\Omega^1(\b{R}_q^N)$ by setting $(\xi^i)^* = \bar\xi^j g_{ji}$.
Then the two calculi are conjugate.

\section{Application of the formalism to $\b R^N_q$}
\label{applic} 

In this section we shall actually construct a frame $\theta^a$ and the 
dual inner derivations $e_a = \ad \lambda_a$ satisfying the conditions 
in Section~\ref{preli2} for $\b{R}^N_q$.

However, before this is possible, we have to enlarge the algebra of $\b R^N_q$.
In Section~\ref{preli2}  we required the center of the algebra $\c{A}$
to be trivial. But the algebra generated by the $x^i$
has a nontrivial center, the element $r^2=g_{kl} x^k x^l$
commutes with all coordinates. For this reason the formalism cannot
be directly applied: with a general Ansatz of the type
\be
\theta^a = \theta^a_i \xi^i                                \label{ansatz}
\ee
it is an immediate result that the condition (\ref{thetaf}) cannot be 
fulfilled for $r^2$.
To solve this problem we add to the algebra an element
$\Lambda$, that we call the ``dilatator'',
and its inverse $\Lambda^{-1}$, 
satisfying the commutation relations
\be
\ba{l}
x^i \Lambda=q \Lambda x^i,                           \\
\xi^i \Lambda=\Lambda \xi^i, \qquad \Lambda d=q d \Lambda.
\nonumber
\ea                                                    \label{xLambda}
\ee
Note that $\Lambda$ does not satisfy the
Leibniz rule $d(fg)=f dg + (df) g \;\forall f,g \in \b{R}^N_q$, so it must 
be interpreted as an element of the $q$-deformed Heisenberg algebra;
in fact  $\Lambda^{-2}$ can be
constructed~\cite{Ogi92} as a simple polynomial in the \uqs-covariant
coordinates and partial derivatives.
The choice in the second line of (\ref{xLambda})
is not unique. Another possibility~\cite{FioMad99}, which was considered
also in \cite{AscCastSca97}, would be $d \Lambda=\Lambda d$ and
$\xi^i\Lambda=q\Lambda \xi^i$, which gives a strict Leibniz rule
also for $f=\Lambda$. However, this is a bit cumbersome because
it implies $df=0$ 
to hold not only for the analogues of the constant functions.

Now, we can proceed with the actual construction of a frame.
We first find one in the larger algebra $\Omega^*(\b R^N_q)
\cocross \uqs$. It is convenient to introduce the
Faddeev-Reshetikin-Takhtadjan generators \cite{FadResTak89} of $\uqs$:
\be
\c{L}^+{}_l^a:=\R^{(1)}\rho_l^a(\R^{(2)})\qquad\qquad
\c{L}^-{}_l^a:=\rho_l^a(\R^{-1}{}^{(1)})\R^{-1}{}^{(2)}, \label{frt}
\ee
where $\uqsp$ and  $\uqsn$ are the Borel subalgebras of \uqs, 
$\R=\R^{(1)} \otimes \R^{(2)} \in \uqsp \otimes \uqsn$ is the 
universal $\R$-matrix for $SO_q(N)$ and $\rho$ the vector 
$N$-dimensional representation of $\uqs$. Then 

\begin{prop} \cite{CerFioMad00}
If we define
\bea
&&\vartheta^a=\Lambda^{-1} \c{L}^-{}^a_l\xi^l \in
\Omega^1(\b R^N_q) \cocross \uqsn        \label{utile1}\\
&&\bar \vartheta^a=\Lambda \c{L}^+{}^a_l \bar\xi^l \in
\bar\Omega^1(\b R^N_q) \cocross \uqsp ,        \label{utile2}
\eea
then  $\vartheta^a$, $\bar \vartheta^a$ constitute (in a generalized
sense) a frame for $\b R^N_q$:
\be
[\vartheta^a,f]=0, \quad [\bar\vartheta^a,f]=0 \quad \forall f \in \b R^N_q.
\ee
Moreover,
\be
[\vartheta^a,\Lambda]=0, \quad [\bar\vartheta^a,\Lambda]=0.
\ee
The commutation relations among the $\vartheta^a$ (resp.
$\bar\vartheta^a$) are  as the ones among the $\xi^i$ (resp. 
$\bar\xi^i$), except for the opposite products:
\be
\begin{array}{ll}
\c{P}_s{}_{ab}^{cd}\vartheta^b\vartheta^a=0,  
\qquad &\c{P}_t{}_{ab}^{cd}\vartheta^b\vartheta^a=0       \cr
 \c{P}_s{}_{ab}^{cd}\bar\vartheta^b\bar\vartheta^a=0,  
\qquad &\c{P}_t{}_{ab}^{cd}\bar\vartheta^b\bar\vartheta^a=0.   
\end{array}                               \label{ththrel}
\ee
\end{prop}

Now, we wish to find a frame in the smaller algebra $\Omega^1(\b R^N_q)$.
To this end, it would suffice to replace in the definitions
(\ref{utile1}-\ref{utile2}) $\c{L}^{\pm}{}^a_l$ by
elements in $\b R^N_q$ which have their same commutation
relations with any element $f\in \b R^N_q$. Such elements exist
since in Ref. \cite{CerFioMad00} we have found homomorphims
$\varphi^{\pm}: \b R^N_q \cocross U_q^{\pm}{\/so(N)} \rightarrow \b R^N_q$
acting as the identity on $\b R^N_q$,
so we just need to replace $\c{L}^{\pm}{}^a_l$ by
$\varphi^{\pm}(\c{L}^{\pm}{}^a_l)$  in order to project 
$\vartheta^a$ and $\bar \vartheta^a$
to frames in $\Omega^1(\b R^N_q)$, 
$\bar\Omega^1(\b R^N_q)$ respectively. Strictly speaking, the homomorphisms
take values  in a slightly enlarged version of the 
algebra of $\b R^N_q$, obtained by adding some new elements
$(r_i)^{\frac 12}$, $i=0,\ldots n$, 
together with their inverses $(r_i)^{-\frac 12}$, 
fulfilling the condition that
$r_i^2=\sum_{k,l=-i}^i g_{kl} x^k x^l$. 
Their commutation relations with the $x^j$ are
are automatically fixed by the latter condition and by the commutation
relations between $r_i^2$ and $x^j$ which can be drawn from 
(\ref{xrel}). They read
\bea
&&x^j r_i=r_i x^j \hbox{ for } |j| \le i, \nn
&&x^j r_i=q r_i x^j \hbox{ for } j < -i, \label{xrrel}\\
&& x^j r_i=q^{-1} r_i x^j \hbox{ for } j > i. \nonumber
\eea
But now for $N$ even the center of the algebra is no longer trivial,
even after the addition of $\Lambda$, because $r_1^{-1} x^{\pm 1}$
commutes with $\Lambda$ as well as with the coordinates.
In this case it is necessary to enlarge the algebra again, and we choose
to add a Drinfeld-Jimbo generator $K=q^{\frac {H_1}2}$ and its
inverse $K^{-1}$, where $H_1$ belongs to the Cartan subalgebra
of $U_qso(N)$ and represents the component of the angular momentum
in the $(-1,1)$-plane.
This new element satisfies the commutation relations
\be
K \Lambda=\Lambda K, \quad K x^{\pm 1}=q^{\pm 1} x^{\pm 1} K, \quad
K x^{\pm i}= x^{\pm i} K \hbox{ for } i>1 \label{xkapparel}
\ee
and we fix its commutation relations with the $1$-forms to be
\be
K \xi^1=q^{\pm 1} \xi^{\pm 1} K, \qquad K \xi^i=\xi^i K \quad
\hbox{ for } i>1.
\ee

\begin{theorem}\cite{CerFioMad00}                        \label{theor1}
A homomorphism $\varphi^-:\b R^N_q \cocross U^-_qso(N)\rightarrow \b R^N_q$ 
can be defined by setting on the generators
\bea
&& \varphi^-(a)=a,         \qquad \forall a\in\b R^N_q, \\[4pt]
&& \varphi^-(\c{L}^-{}^i_j)=g^{ih}\Lambda^{-1}[\lambda_h,x^k]g_{kj},
\eea
with
\be
\begin{array}{ll}
\lambda_0=\gamma_0 \Lambda (x^0)^{-1}
&\quad\mbox{for $N$ odd,} \\[6pt]
\lambda_{\pm 1}=\gamma_{\pm 1} \Lambda (x^{\pm 1})^{-1} K^{\mp 1} 
&\quad\mbox{for $N$ even,} \\[6pt]
\lambda_a=\gamma_a \Lambda r_{|a|}^{-1}r_{|a|-1}^{-1} x^{-a}
&\quad\mbox{otherwise,} 
\end{array}                                             \label{deflambda}
\ee
and $\gamma_a \in \b{C}$ normalization constants fulfilling
the conditions
\be
\begin{array}{ll}
\gamma_0 = -q^{-\frac{1}{2}} h^{-1} &\quad\mbox{for $N$ odd,} \\[6pt]
\gamma_1 \gamma_{-1}=
\left\{\begin{array}{l}
-q^{-1} h^{-2}\\
k^{-2}
\end{array}\right.
&\quad\!\begin{array}{l}
\mbox{for $N$ odd,} \\
\mbox{for $N$ even,}
\end{array}\\[8pt]
\gamma_a \gamma_{-a} =
-q^{-1} k^{-2} \omega_a \omega_{a-1} &\quad\mbox{for $a>1$}. \nonumber
\end{array}                                               \label{gamma}
\ee
\end{theorem}
Here and in the sequel we set 
\be
h:=\sqrt{q}-1/\sqrt{q},\qquad k:=q-q^{-1},\qquad 
\omega_i := q^{\rho_i}+q^{-\rho_i}.
\ee
The relation (\ref{gamma}) fixes only the product $\gamma_a \gamma_{-a}$,  
$\gamma_0^2$ for $N$ odd and $\gamma_1 \gamma_{-1}$ for $N$ even
are positive real numbers, while all the remaining products 
$\gamma_a \gamma_{-a}$ are negative.

\begin{theorem}\cite{CerFioMad00}                   \label{theor2}
A homomorphism $\varphi^+:\b R^N_q \cocross \uqsp \rightarrow \b R^N_q$ by
setting on the generators
\bea
&& \varphi^+(a)=a,           \qquad \forall a\in\b R^N_q, \\[4pt]
&& \varphi^+(\c{L}^+{}^i_j)=g^{ih}\Lambda[\bar\lambda_h,x^k]g_{kj},
\eea
with
\be
\begin{array}{ll}
\bar\lambda_0=\bar\gamma_0 \Lambda^{-1} (x^0)^{-1}
&\quad\mbox{for $N$ odd,} \\[6pt]
\bar\lambda_{\pm 1} = 
\bar\gamma_{\pm 1} \Lambda^{-1} (x^{\pm 1})^{-1} K^{\pm 1} 
&\quad\mbox{for $N$ even,} \\[6pt]
\bar\lambda_a = 
\bar\gamma_a \Lambda^{-1} r_{|a|}^{-1}r_{|a|-1}^{-1} x^{-a}
&\quad\mbox{otherwise,} 
\end{array}                                         \label{defbarlambda}
\ee
and $\bar\gamma_a \in \b{C}$ normalization constants fulfilling
the conditions
\be
\begin{array}{ll}
\bar \gamma_0 = q^{\frac{1}{2}} h^{-1}
&\quad\mbox{for $N$ odd,} \\[6pt]
\bar \gamma_1 \bar \gamma_{-1} =
\left\{
\begin{array}{l}
-q h^{-2}\\
k^{-2}
\end{array}
\right.
&\quad\!\begin{array}{l}
\mbox{for $N$ odd,} \\
\mbox{for $N$ even,}
\end{array}\\[8pt]
\bar \gamma_a \bar \gamma_{-a}=
-q k^{-2} \omega_a \omega_{a-1}
&\quad\mbox{for $a>1$}.         \nonumber
\end{array}                                             \label{bargamma}
\ee
\end{theorem}
Under which circumstances can these homomorphisms be `glued' together to
a homomorphism of the whole of \uqs ? The answer is given by 

\begin{theorem}\cite{CerFioMad00}
In the case of odd $N$ we one can define a homomorphism
$\varphi: \uqs \cross\b R^N_q\rightarrow \b R^N_q$
by setting on the generators
\bea
&& \varphi(a)=a\qquad \forall a\in\b R^N_q \\[4pt]
&& \varphi(\c{L}^-{}^i_j)=g^{ih}\Lambda^{-1}[\lambda_h,x^k]g_{kj},\\[4pt]
&& \varphi(\c{L}^+{}^i_j)=g^{ih}\Lambda[\bar\lambda_h,x^k]g_{kj},
\eea
with $\lambda_j,\bar\lambda_j$ defined as in (\ref{deflambda}),
(\ref{defbarlambda}) and with coefficients given by
\be
\begin{array}{lr}
\gamma_0 = -q^{-\frac 12}h^{-1}\\[4pt]
\gamma_1^2 = -q^{-2} h^{-2} \\[4pt]
\gamma_a^2 = -q^{-2}\omega_a \omega_{a-1} k^{-2} 
&\quad \mbox{for $a>1$}\\[4pt]
\gamma_a = q \gamma_{-a}  
& \quad \mbox{for $a \le 1$} \\[4pt]
\bar \gamma_a=-q \gamma_a.
\end{array}                                              \label{fixgamma} 
\ee
\end{theorem}
Notice that the $\gamma_a,\bar\gamma_a$ for $a \neq 0$ are imaginary
and fixed only up to a sign. Therefore, the
homomorphism $\varphi$ does not preserve the star structure of 
the real section $U_qso(N,\b{R})$
for $q\in\b{R}^+$.

In the case of even $N$ it is not possible to extend 
$\varphi^{\pm}$ to a homomorphism $\varphi$ from the whole of \uqs,
due to the necessity of adding $K^{\pm 1}$ to the algebra.
In other words in this case there is one of the generators of the Cartan
subalgebra of \uqs which cannot be realized in $\b R^N_q$ and this has a 
consequence that it is not possible to satisfy the would-be
$\varphi$-image of the commutation relations between $\c{L}^+{}i_j$
and $\c{L}^-{}h_k$, so that the homomorphisms cannot be extended to 
the whole \uqs.

The fact that the homomorphism
$\varphi^{\pm}$ or $\varphi$ exist can be interpreted as the
existence of a local realization of $U_q^{\pm}so(N)$ or \uqs.

Summing up, for the frames of the two calculi we explicitly find:
\be
\theta^a=\theta^a_l\xi^l=
\Lambda^{-2}g^{ab}[\lambda_b,x^j]g_{jl}\xi^l, \qquad
\bar\theta^a=\bar\theta^a_l\bar\xi^l=
\Lambda^2g^{ab}[\bar\lambda_b,x^j]g_{jl}\bar\xi^l.
\ee
The $\theta^a$ commute both with the coordinates and 
with $\Lambda$. 

{From} the preceding theorems one derives the 
commutation relations 
\begin{theorem}\cite{CerFioMad00}
\bea
\c{P}_a{}^{ab}_{cd} \lambda_a \lambda_b=0, \qquad\qquad
\c{P}_a{}^{ab}_{cd} \bar\lambda_a \bar \lambda_b=0 \\
{\c{P}_a}^{ab}_{cd} \theta^c \theta^d=\theta^a \theta^b, \qquad
\c{P}_a{}^{ab}_{cd} \bar\theta^c \bar\theta^d = \bar\theta^a \bar\theta^b.
\eea
\end{theorem}
Hence the matrix $P$ of equations (\ref{prod}) and (\ref{manca})
is the $q$-deformed antisymmetric projector $\c{P}_a$.
The elements $\lambda_a$, $\bar \lambda_a$ satisfy the same commutation
relation as the $x^i$, while the $\theta^a, \bar \theta^a$ those
satisfied by the $\xi^i$.

It is a very interesting result that applying a condition like
(\ref{thetaf}), which does not depend on the quantum group covariance 
of $\b R^N_q$,  to the quantum group covariant differential calculi,
the components $\theta^a_i$ of the frame automatically carry so to
say a `local realization' of \uqs.

The Dirac operator~\cite{Con94} of \ref{dirac} is easily verified to be
given by~\cite{Ste96,Zum97}.
\bea
\theta&= \omega_n q^{\frac{N}{2}} k^{-1} r^{-2}
g_{ij} x^i \xi^j, &\mbox{for } \Omega^1(\b{R}_q^N),\\
\bar\theta &= -\omega_n q^{-\frac{N}{2}} k^{-1} r^{-2} g_{ij} x^i \bar\xi^j
&\mbox{for } \bar\Omega^1(\b{R}_q^N).
\eea

Finally we summarize the results found in \cite{CerFioMad00} for the metric
and the covariant derivative for each of the two calculi $\Omega(\b{R}^N_q)$
and $\Omega^*(\b{R}^N_q)$. In the $\theta^a$, $\bar \theta^a$ basis 
respectively the actions of $g$ and $\sigma$ are
\bea
\sigma (\theta^a\otimes \theta^b) =S^{ab}{}_{cd}\,\theta^c\otimes \theta^d,
&\:\: g(\theta^a\otimes \theta^b) = g^{ab} &\:\:
\mbox{for }\Omega^*(\b{R}^N_q),\\
\sigma (\bar\theta^a\otimes \bar\theta^b) =
\bar S^{ab}{}_{cd}\,\bar\theta^c\otimes \bar\theta^d,
&\:\:g(\bar\theta^a\otimes\bar\theta^b) = g^{ab}
&\:\:\mbox{for }\bar\Omega^*(\b{R}^N_q).
\eea
while in the $\xi^i$, $\bar\xi^i$ basis we get
\bea
g(\xi^i \otimes \xi^j)= g^{ij} \Lambda^2 &
\sigma(\xi^i \otimes \xi^j)=S^{ij}_{hk} \xi^h \otimes \xi^k &
\mbox{for }\Omega^*(\b{R}^N_q),\\
g(\bar \xi^i\otimes \bar\xi^j) = g^{ij} \Lambda^{-2},
&\sigma (\bar \xi^i\otimes \bar \xi^j) =
\bar S^{ij}{}_{hk}\,\bar\xi^h\otimes \bar\xi^k.
&\mbox{for }\bar\Omega^*(\b{R}^N_q).
\eea
Unfortunately, it is not possible to satisfy simultaneously the metric
compatibility condition (\ref{met-comp}) and the bilinearity condition
for the torsion (\ref{torsion}). The best we can do is
to weaken the compatibility condition to a condition of proportionality.
Then for each calculus we find the two solutions for $\sigma$:
\bea
S = q\hat R, & S = (q\hat R)^{-1} &\mbox{for }\Omega^*(\b{R}^N_q),\\
\bar S = q\hat R, & \bar S = (q\hat R)^{-1}
&\mbox{for }\bar\Omega^*(\b{R}^N_q).
\eea
This implies that the covariant derivatives and metric are compatible
only up to a conformal factor:
\be
S^{ae}_{df} g^{fg} S^{cb}_{eg}=q^{\pm 2} g^{ac} \delta^b_d, \qquad
\bar S {}^{ae}{}_{df} g^{fg} \bar S {}^{cb}{}_{eg} 
= q^{\pm 2} g^{ac} \delta^b_d.
\ee
The linear curvatures associated to the covariant derivatives 
defined by (\ref{covdev0}) for each the calculi vanish, consistently
with the fact that $\b R^N_q$ should be flat.


\begin{thebibliography}{10}

\bibitem{CerFioMad00}
B. L. Cerchiai, G. Fiore, J. Madore, ``Geometrical Tools
for Quantum Euclidean Spaces'', preprint 99-52, Dip. Mat. e 
Appl., Univ. di Napoli, LMU-TPW 99-17, MPI-PhT/99-45, 
math.QA/0002007.

\bibitem{Sny47}
H.~Snyder, 
{\em Phys.\ Rev.} {\bf 71} (1947) 38;
{\em Phys.\  Rev.} {\bf 72} (1947) 68.

\bibitem{DimMad96}
A.~Dimakis, J.~Madore, 
{\em  J.~Math.\ Phys.} {\bf 37} (1996), no.~9, 4647.

\bibitem{FadResTak89}
L.~Faddeev, N.~Reshetikhin, L.~Takhtadjan, 
{\em Leningrad Math. J.} {\bf 1} (1990) 193.

\bibitem{Mad95}
J.~Madore, ``An Introduction to Noncommutative Differential
Geometry and its Physical Applications'',
London Mathematical Society Lecture Note Series {\bf 257},
Cambridge University Press, 2nd Ed., 1999


\bibitem{CarSchWat91}
U.~Carow-Watamura, M.~Schlieker, S.~Watamura, 
{\em Z.~Physik~C - Particles and Fields} {\bf 49} (1991) 439.

\bibitem{Ogi92}
O.~Ogievetsky, 
{\em Lett.\ Math.\ Phys.} {\bf 24} (1992) 245.

\bibitem{WesZum90}
J.~Wess, B.~Zumino, 
{\em Nucl.\ Phys.\ (Proc.\ Suppl.)} {\bf 18B} (1990) 302.

\bibitem{Con94}
A.~Connes, ``Noncommutative Geometry'',
\newblock Academic Press, 1994.

\bibitem{Mou95}
J.~Mourad,
{\em Class. Quant. Grav.} {\bf 12} (1995), 965

\bibitem{Zum97}
C.~S.~Chu, P. M.~Ho, B.~Zumino,
``Some complex quantum manifolds and their geometry'', in 
{\it Quantum Fields and Quantum Space Time},
G.'t Hooft, A.Jaffe, G.Mack, P. Mitter, R. Stora (Eds.),
Plenum Press, NY, NATO ASI Series {\bf 364} (1997) 283. 

\bibitem{Ste96}
H.~Steinacker, 
{\em J.~Math.\ Phys.} {\bf 37} (1996) 7438.

\bibitem{Wor}
S.~Woronowicz, 
{\em Commun.\ Math.\ Phys.}  {\bf 111} (1987) 613; 
{\em Commun.\ Math.\ Phys.} {\bf 122} (1989) 125.  

\bibitem{AscCast96} P. Aschieri, L. Castellani,
{\em  Int. J. Mod. Phys.} {\bf A11} (1996), 4513.

\bibitem{AscCastSca97} P. Aschieri, L. Castellani, A. M. Scarfone
{\em  Eur. Phys. J.} {\bf C7} (1999), 159-175.


\bibitem{FioMad00}
G. Fiore, J. Madore, ``Geometrical Issues for the
3-dim Quantum Euclidean Space'', preprint 99-53, Dip. Mat. e Appl. Napoli, 
math.QA/00029191, talk given at the VI. Wigner Symposium,
Istanbul, Turkey, August 1999, to appear in the proceedings

\bibitem{GroMadSte00}
H. Grosse, J. Madore, H. Steinacker, ``Field Theory on the $q$-deformed
fuzzy sphere I'', preprint LMU-TPW 00-13, UWThPh-2000-19, LPTE Orsay 00/54,
hep-th/0005273, 2000

\bibitem{FioMad98}
G.~Fiore, J.~Madore, ``Leibniz rules and reality conditions'', preprint 
98-13, Dip. Mat. e Appl., Univ. di Napoli, math.QA/9806071.

\bibitem{DubMadMasMou96}
M.~Dubois-Violette, J.~Madore, T.~Masson, J.~Mourad, 
{\em J.~Math.\ Phys.} {\bf 37} (1996), 4089.

\bibitem{DubMadMasMou95}
M.~Dubois-Violette, J.~Madore, T.~Masson, J.~Mourad, 
{\em Lett.\ Math.\ Phys.} {\bf 35} (1995), 351.

\bibitem{FioMad99}
G.~Fiore, J.~Madore, 
{\em J. Geom. Phys.} {\bf 33} (2000), 257-287.

\end{thebibliography}
\end{document}